# Decentralized Static Output Feedback Controller Design for Large Scale Switched T-S Systems

Dalel Jabri[1], Djamel Eddine Chouaib Belkhiat[1], Kevin Guelton[2], Noureddine Manamanni[2]

[1] Laboratoire DAC HR
Ferhat Abbas University, Setif1
Campus El Bez, Sétif 19000, Algeria
djamel.belkhiat@univ-setif.dz, dalel.jabri@yahoo.fr

[2]University of Reims Champagne-Ardenne
CReSTIC EA3804, Moulin de la Housse
BP 1039, 51687 Reims Cedex 2, France.
Kevin.guelton@univ-reims.fr, Noureddine.Manamanni@univ-reims.fr

*Abstract*— **This paper investigates the design of decentralized output-feedback controllers for a class of a large scale switched nonlinear systems under arbitrary switching laws. A global large scale switched system can be split into a set of smaller interconnected switched Takagi Sugeno fuzzy subsystems. Then, in order to stabilize the overall closed-loop system, a set of switched non-PDC static output controllers is employed. The latter is designed based on Linear Matrix Inequality (LMI) conditions obtained from a multiple switched non quadratic-like Lyapunov candidate function. The controllers proposed herein are synthesized to satisfy $H_\infty$ performance for disturbance attenuation. Finally, a numerical example is proposed to illustrate the effectiveness of the suggested decentralized switched controller design approach.**

## I. INTRODUCTION

During the latter decades, several complex systems are appeared to meet the specific needs of the world population. In this context, we can quote as examples networked power systems, water transportation networks, traffic systems, as well as other systems in various fields. Generally speaking, establishing mathematical models for these systems is a complex task, especially when the system is considered as a whole. Hence, to overcome such difficulties, an alternative to global modelling approach can be considered. It consists in splitting the overall large-scale system in a finite set of interconnected low-order subsystems [1].

Among these complex systems, switched interconnected large-scale system have attracted considerable attention since they provide a convenient modelling approach for many physical systems that can exhibit both continuous and discrete dynamic behaviour. In this context, several studies dealing with the stability analysis and stabilization issues for both linear and nonlinear switched interconnected large-scale systems have been explored [1]-[7]. Hence, the main challenge to deal with such problems consists in determining the conditions ensuring the stability of the whole system with consideration to the interconnections effects between its subsystems. Nevertheless, few works based on the approximation property of Takagi-Sugeno (TS) fuzzy models for nonlinear problems, have been achieved to deal with the stabilization of continuous-time large-scale switched nonlinear systems [3], [7], [8].

Hence, this paper presents the design of decentralized static output feedback controllers for a class of switched Takagi-Sugeno interconnected large-scale system with external bounded disturbances. More specifically, the primary contribution of this paper consists in proposing a LMI-based methodology, in the non quadratic framework, for the design of robust decentralized switched non-PDC controllers for a class of large scale switched nonlinear systems under arbitrary switching laws.

The remainder of the paper is organized as follows. Section 2 presents the considered class of switched Takagi-Sugeno interconnected large-scale system, followed by the problem statement. The design of the considered decentralized and switched static output feedback non-PDC controllers is presented in section 3. A numerical example is proposed to illustrate the efficiency of the proposed approach in section 4. The paper ends with conclusions and references.

## II. PROBLEM STATEMENT AND PRELIMINARIES

Let us consider the class of nonlinear hybrid systems $S$ composed of $n$ continuous time switched nonlinear subsystems $S_i$ represented by switched TS models. The $n$ state equations of the whole interconnected switched fuzzy system $S$ are given as follows; for $i=1,...,n$:

$$\dot{x}_i(t) = \sum_{j_i=1}^{m_i}\sum_{s_{j_i}=1}^{r_{j_i}} \xi_{j_i}(t) h_{s_{j_i}}(z_{j_i}(t)) \begin{bmatrix} A_{s_{j_i}} x_i(t) + B_{s_{j_i}} u_i(t) + B^w_{s_{j_i}} w_i(t) \\ + \sum_{\alpha=1,\alpha\neq i}^{n} \left( F_{i,\alpha,s_{j_i}} x_\alpha(t) + B^{w_\alpha}_{s_{j_i}} w_\alpha(t) \right) \end{bmatrix} \quad (1)$$

$$y_i(t) = \sum_{j_i=1}^{m_i}\sum_{s_{j_i}=1}^{q_{j_i}} \xi_{j_i}(t) h_{s_{j_i}}(z_{j_i}(t)) C_{s_{j_i}} x_i(t) \quad (2)$$

where $x_i(t) \in \mathbb{R}^{\eta_i}$, $y_i(t) \in \mathbb{R}^{\rho_i}$, $u_i(t) \in \mathbb{R}^{\upsilon_i}$ represent respectively the state, the measurement (output) and the input vectors associated to the $i^{th}$ subsystem. $w_i(t) \in \mathbb{R}^{\upsilon_i}$ is an $L_2$-norm-bounded external disturbance associated to the $i^{th}$ subsystem. $m_i$ is the number of switching modes of the $i^{th}$ subsystem. $r_{j_i}$ is the number of fuzzy rules associated to the $i^{th}$ subsystem in the $j_i^{th}$ mode; for $i=1,\ldots,n$, $j_i = 1,\ldots,m_i$ and $s_{j_i} = 1,\ldots,r_{j_i}$, $A_{s_{j_i}} \in \mathbb{R}^{\eta_i \times \eta_i}$, $B_{s_{j_i}} \in \mathbb{R}^{\eta_i \times \upsilon_i}$, $B_{s_{j_i}}^w \in \mathbb{R}^{\eta_i \times \upsilon_i}$, $C_{l_{j_i}} \in \mathbb{R}^{\rho_i \times \eta_i}$, $F_{i,\alpha,s_{j_i}} \in \mathbb{R}^{\eta_i \times \eta_\alpha}$ and $B_{s_{j_i}}^{w_\alpha} \in \mathbb{R}^{\eta_i \times \upsilon_\alpha}$ are constant matrices describing the local dynamics of each polytops; $F_{i,\alpha,s_{j_i}}$, $B_{s_{j_i}}^{w_\alpha}$ express the interconnections between subsystems. $z_{j_i}(t)$ are the premises variables and $h_{s_{j_i}}(z_{j_i}(t))$ are positive membership functions satisfying the convex sum proprieties $\sum_{s_{j_i}=1}^{r_{j_i}} h_{s_{j_i}}(z_{j_i}(t)) = 1$; $\xi_{j_i}(t)$ is the switching rules of the $i^{th}$ subsystem, considered arbitrary but assumed to be real time available, these are defined such that the active system in the $l_i^{th}$ mode lead to:

$$\begin{cases} \xi_{j_i}(t) = 1 & \text{if } j_i = l_i \\ \xi_{j_i}(t) = 0 & \text{if } j_i \neq l_i \end{cases} \quad (3)$$

**Notations**: In order to lighten the mathematical expression, one assumes the scalar $\underline{N} = \frac{1}{n-1}$, the index $j_i$ denote the switched modes associated to the $i^{th}$ subsystem. The premises entries $z_{j_i}$ will be omitted when there is no ambiguities and the following notations will be employed in the sequel:

$G_{hj} = \sum_{s_{j_i}=1}^{r_{j_i}} h_{s_{j_i}} G_{s_{j_i}}$ and $Y_{hj,hj} = \sum_{s_{j_i}=1}^{r_{j_i}} \sum_{k_{j_i}=1}^{r_{j_i}} h_{s_{j_i}} h_{k_{j_i}} Y_{s_{j_i},k_{j_i}}$. For matrices of appropriate dimensions we will denote: $\dot{X}_{hj} = \frac{dX_{h_{j_i}}}{dt}$ and $\left(\dot{X}_{hj}\right)^{-1} = \frac{d\left(X_{h_{j_i}}\right)^{-1}}{dt}$. As usual, a star (*) indicates a transpose quantity in a symmetric matrix. The time $t$ will be omitted when there is no ambiguity. However, one denotes $t_{j \to j^+}$ the switching instants of the $i^{th}$ subsystem between the current mode $j$ (at time $t$) and the upcoming mode $j^+$ (at time $t^+$), therefore:

$$\begin{cases} \xi_j(t) = 1 \\ \xi_{j^+}(t) = 0 \end{cases} \text{ and } \begin{cases} \xi_j(t^+) = 0 \\ \xi_{j^+}(t^+) = 1 \end{cases} \quad (4)$$

In the sequel, we will deal with the design of static output-feedback controllers with disturbance attenuation for the considered class of large-scale system $S$. For that purpose, a set of decentralized static output feedback switched non-PDC control laws is proposed as; for $i=1,\ldots,n$:

$$u_i(t) = \sum_{j_i=1}^{m_i} \sum_{k_{j_i}=1}^{r_{j_i}} \xi_{j_i}(t) h_{s_{j_i}}(z_{j_i}(t)) K_{k_{j_i}} \left(\sum_{s_{j_i}=1}^{r_{j_i}} h_{s_{j_i}}(z_{j_i}(t)) X_{s_{j_i}}^9 \right)^{-1} y_i(t) \quad (5)$$

where $K_{k_{j_i}}$ and $X_{s_{j_i}}^9 = \left(X_{s_{j_i}}^9\right)^T > 0$ are the non-PDC gain matrices to be synthesized.

Hence, substituting (5) into (1), one expresses the overall closed-loop dynamics $S_{cl}$ described by; for $i=1,\ldots,n$:

$$\dot{x}_i = \sum_{j=1}^{m_i} \xi_j \left\{ \begin{array}{l} \left[A_{hj} + B_{hj} K_{hj} \left(X_{hj}^9\right)^{-1} C_{hj}\right] x_i + B_{s_{j_i}}^w w_i(t) \\ + \sum_{\alpha=1,\alpha \neq i}^n \left(F_{i,\alpha,hj} x_\alpha + B_{s_{j_i}}^{w_\alpha} w_\alpha(t)\right) \end{array} \right\} \quad (6)$$

Thus, the problem considered in this study can be resumed as follows:

**Problem 1**: The objective is to design the static output feedback controllers (5) such that the switched TS interconnected large-scale system (1)-(2) rises a closed-loop robust $H_\infty$ output-feedback stabilization performance.

**Definition 1**: The switched interconnected large-scale system (1)-(2) is said to have a robust $H_\infty$ output-feedback performance if the following conditions are satisfied:

- **Condition 1** *(Stability condition)*: With zero disturbances input condition, i.e. $w_i \equiv 0$ for $i=1,\ldots,n$, the closed-loop dynamics (6) is stable.

- **Condition 2** *(Robustness condition)*: For all non-zero $w_i \in L_2[0 \; \infty)$, under zero initial condition $x_i(t_0) \equiv 0$, the following $H_\infty$ criterion holds for $i=1,\ldots,n$,

$$J_i = \int_0^{+\infty} x_i^T x_i dt \leq \varsigma_i^2 \int_0^{+\infty} \left(w_i^T w_i + \sum_{\alpha=1,\alpha \neq i}^n w_\alpha^T w_\alpha \right) dt \quad (7)$$

where $\varsigma_i^2$ is positive scalars which represents the disturbance attenuation level associated to the $i^{th}$ subsystem.

From the closed-loop dynamics (6), it can be seen that several crossing terms among the gain controllers $K_{hj}$ and the system's matrices $\left(B_{hj} K_{hj} \left(X_{hj}^9\right)^{-1} C_{hj}\right)$ are present. Hence, in view of the wealth of interconnections characterizing our system, these crossing terms lead surely to very conservative conditions for the design of the proposed controller. In order to decouple the crossing terms $\left(B_{hj} K_{hj} \left(X_{hj}^9\right)^{-1} C_{hj}\right)$ appearing in the equation (6) and to provide LMI-based design conditions, an interesting property called 'descriptor redundancy' can be considered [10]. Thus, the closed-loop dynamics (6) can be alternatively expressed as follows. First, from (2), we introduce null terms and it yields, for $i=1,\ldots,N$:

$$0\dot{y}_i = -y_i + C_{hj} x_i, \quad (8)$$

$$0 = u_i - K_{hj} \left(X_{hj}^9\right)^{-1} y_i \quad (9)$$

Then, by considering the following augmented variables $\tilde{x}_i^T = \begin{bmatrix} x_i^T & y_i^T & u_i^T \end{bmatrix}$, $\tilde{x}_\alpha^T = \begin{bmatrix} x_\alpha^T & y_\alpha^T & u_\alpha^T \end{bmatrix}$, $\tilde{w}_{i,\alpha}^T = \begin{bmatrix} w_i^T & w_\alpha^T \end{bmatrix}$, the large-scale system (1)-(9) and the TS controllers (5) can be combined as follows to express the closed-loop dynamics.
For $i = 1, \ldots, N$:

$$E\dot{\tilde{x}}_i = \tilde{A}_{hj,hj}\tilde{x}_i + \sum_{\alpha=1,\alpha \neq i}^{n} \left( \tilde{F}_{i,\alpha,hj}\tilde{x}_\alpha + \tilde{B}_{hj}^{w\alpha}\tilde{w}_{i,\alpha} \right) \tag{10}$$

with $E = \begin{bmatrix} I & 0 & 0 \\ 0 & 0 & 0 \\ 0 & 0 & 0 \end{bmatrix}$, $\tilde{A}_{hj,hj} = \begin{bmatrix} A_{hj} & 0 & B_{hj} \\ 0 & -I & K_{hj}\left(X_{hj}^9\right)^{-1} \\ C_{hj} & 0 & -I \end{bmatrix}$

$\tilde{F}_{i,\alpha,hj} = \begin{bmatrix} F_{i,\alpha,hj} & 0 & 0 \\ 0 & 0 & 0 \\ 0 & 0 & 0 \end{bmatrix}$ and $\tilde{B}_{hj}^{w\alpha} = \begin{bmatrix} NB_{hj}^w & 0 \\ 0 & B_{hj}^{w\alpha} \\ 0 & 0 \end{bmatrix}$.

Note that the system (10) is a large scale switched descriptor. Hence, it is worth pointing out that the output-feedback stabilization problem of the system (1)-(2) can be converted into the stabilization problem of the augmented system (10).

*Remark*: If it may be difficult to work on the first formulation of the closed-loop dynamics (6) due to the large number of crossing terms, the goal of our study can now be achieved by considering the augmented closed-loop dynamics (10) expressed in the descriptor form. In this context, the second condition of the definition 1, given by equation (7), can be reformulated as follows:

$$\int_0^{+\infty} \tilde{y}_i^T Q \tilde{y}_i dt \leq \varsigma_i^2 \int_0^{+\infty} \sum_{\alpha=1,\alpha \neq i}^{n} \tilde{w}_{i,\alpha}^T \Xi \tilde{w}_{i,\alpha} dt \tag{11}$$

with $\Xi = \begin{bmatrix} NI & 0 \\ 0 & I \end{bmatrix}$ $Q = \begin{bmatrix} 0 & 0 & 0 \\ 0 & I & 0 \\ 0 & 0 & 0 \end{bmatrix}$

To conclude this preliminary section, let us introduce the following lemma which will be useful in the main result demonstration.

*Lemma* [9]: Let us consider two matrices $A$ and $B$ with appropriate dimensions and a positive scalar $\tau$, the following inequality is always satisfied:

$$A^T B + B^T A \leq \tau A^T A + \tau^{-1} B^T B \tag{12}$$

### III. LMI BASED DECENTRALIZED CONTROLLER DESIGN

In this section, our results on the design of static output feedback $H_\infty$ decentralized switched non-PDC controllers (5), which ensures the closed-loop stability of (6) and the $H_\infty$ disturbance attenuation performance (11) are presented. The main result is summarized in the following theorem.

***Theorem*** : Assume that for each subsystem $i$ of (1), the active mode is denoted by $j_i$ and, for $j_i = 1, \ldots, m_i$ and $s_{j_i} = 1, \ldots, r_{j_i}$, $\dot{h}_{s_{j_i}}(z(t)) \geq \lambda_{s_{j_i}}$. The overall interconnected switched Takagi-Sugeno system (1)-(2) is stabilized by a set of $n$ decentralized static output feedback switched non-PDC control laws (5) according to the definition 1, if there exists, for all combinations of $i = 1, \ldots, n$, $j_i = 1, \ldots, m_i$, $j_i^+ = 1, \ldots, m_i$, $s_{j_i} = 1, \ldots, r_{j_i}$, $k_{j_i} = 1, \ldots, r_{j_i}$, $k_{j_i}^1 = 1, \ldots, r_{j_i}$ and $l_{j_i} = 1, \ldots, r_{j_i}$, the matrices $X_{k_{j_i}}^1 = \left(X_{k_{j_i}}^1\right)^T > 0$, $X_{k_{j_i}}^5 = \left(X_{k_{j_i}}^5\right)^T > 0$; $X_{k_{j_i}}^9 = \left(X_{k_{j_i}}^9\right)^T > 0$ $W_{s_{j_i}s_{j_i}k_{j_i}}^1$, $K_{k_{j_i}}$ and the scalars, $\tau_{1,i}$, ... $\tau_{i-1,i}$, $\tau_{i+1,i}$,..., $\tau_{n,i}$ (excepted $\tau_{i,i}$ which don't exist since there is no interaction between a subsystem and himself), such that the LMIs described by (11), (13) and (14) are satisfied.

$$X_{k'_{j_i}}^1 + W_{s_{j_i}k_{j_i}l_{j_i}} > 0 \tag{13}$$

$$\begin{bmatrix} -\mu_{j_i \to j_i^+} X_{k_{j_i}}^1 & X_{k_{j_i}}^1 \\ X_{k_{j_i}}^1 & -X_{k_{j_i^+}}^1 \end{bmatrix} \leq 0 \tag{14}$$

$$\left( \begin{array}{c|cc|c} \Lambda_{s_{j_i}k_{j_i}} & & & * \\ \hline X_{k_{j_i}} & -\tau_{\alpha,i}I & & \\ \hline \left(\tilde{B}_{s_{j_i}}^{w,\alpha}\right)^T & 0 & -N\varsigma_i^2 I & 0 \\ & & 0 & -I \end{array} \right) < 0 \tag{15}$$

$$\left( \begin{array}{c|ccc|ccc} \Gamma_{s_{j_i}l_{j_i}k_{j_i}} & X_{k_{j_i}} & \cdots & \cdots & X_{k_{j_i}} & X_{k_{j_i}} & \cdots & \cdots & X_{k_{j_i}} \\ \hline X_{k_{j_i}} & -\tau_{1,i}I & 0 & \cdots & 0 & 0 & \cdots & \cdots & 0 \\ \vdots & 0 & \ddots & \ddots & \vdots & \vdots & \ddots & 0 & \vdots \\ \vdots & \vdots & \ddots & \ddots & 0 & \vdots & 0 & \ddots & \vdots \\ X_{k_{j_i}} & 0 & \cdots & 0 & -\tau_{i-1,i}I & 0 & \cdots & \cdots & 0 \\ \hline X_{k_{j_i}} & 0 & \cdots & 0 & -\tau_{i+1,i}I & 0 & \cdots & 0 \\ \vdots & \vdots & \ddots & 0 & \vdots & 0 & \ddots & \ddots & \vdots \\ \vdots & \vdots & 0 & \ddots & \vdots & \vdots & \ddots & \ddots & 0 \\ X_{k_{j_i}} & 0 & \cdots & \cdots & 0 & 0 & \cdots & 0 & -\tau_{n,i}I \end{array} \right) < 0 \tag{16}$$

With $\Lambda_{s_{j_i}l_{j_i}k_{j_i}} = \begin{bmatrix} \Gamma_{s_{j_i}l_{j_i}k_{j_i}} & & * \\ & & * \\ & & * \\ \hline X_{k_{j_i}}^2 & 0 & 0 & -I \end{bmatrix}$, $\tilde{B}_{kj}^{w\alpha} = \begin{bmatrix} NB_{kj}^w & 0 \\ 0 & B_{kj}^{w\alpha} \\ 0 & 0 \end{bmatrix}$,

$\Phi_{s_{j_i}l_{j_i}k_{j_i}k'_{j_i}} = \sum_{l_{j_i}=1}^{r_{j_i}} \lambda_{l_{j_i}} \left( X_{l_{j_i}}^1 + W_{s_{j_i}k_{j_i}k'_{j_i}} \right)$, $X_{k_{j_i}} = \begin{bmatrix} X_{k_{j_i}}^1 & 0 & 0 \\ 0 & X_{k_{j_i}}^5 & 0 \\ 0 & 0 & X_{k_{j_i}}^9 \end{bmatrix}$ and

$$\Gamma_{s_{j_i} l_{j_i} k_{j_i}} = \begin{bmatrix} X^1_{k_{j_i}} A_{s_{j_i}}^T + A_{s_{j_i}} X^1_{k_{j_i}} \\ + \tau_{i,\alpha} F_{i,\alpha,s_{j_i}} F_{i,\alpha,s_{j_i}}^T - \Phi_{s_{j_i} l_{j_i} k_{j_i} k'_{j_i}} & (*) & (*) \\ 0 & -X^5_{k_{j_i}} - X^5_{k_{j_i}} & (*) \\ X^9_{k_{j_i}} (B_{s_{j_i}})^T + C_{s_{j_i}} X^1_{k_{j_i}} & (K_{l_{j_i}})^T & -X^9_{k_{j_i}} - X^9_{k_{j_i}} \end{bmatrix}$$

**Proof**: Indeed, the present proof is divided in two parts corresponding to the condition 1 and 2 given in the definition 1.

***Condition 1** (Stability condition)*: With zero disturbances input condition $\tilde{w}_{i,\alpha} \equiv 0$, for $i = 1, \ldots, N$. Let us define the following multiple switched non-quadratic Lyapunov candidate function:

$$V(x_1, x_2, \ldots, x_n) = \sum_{i=1}^{n} \sum_{j_i=1}^{m_i} \xi_{j_i} v_{j_i}(x_i) > 0 \quad (17)$$

where $v_{j_i} = \tilde{x}_i^T E(X_{hj})^{-1} \tilde{x}_i = \tilde{x}_i^T E \left( \sum_{s_{j_i}=1}^{r_{j_i}} h_{s_{j_i}} X_{s_{j_i}} \right)^{-1} \tilde{x}_i$ and with

$EX_{hj} = X_{hj} E > 0$, $X^1_{hj} = X^{1T}_{hj}$, $X_{hj} = \begin{bmatrix} X^1_{hj} & 0 & 0 \\ 0 & X^5_{hj} & 0 \\ 0 & 0 & X^9_{hj} \end{bmatrix}$.

The augmented system (10), and implicitly the closed-loop interconnected switched system (6), is asymptotically stable if:

$$\forall t \neq t_{j \to j^+}, \quad \dot{V}(x_1, x_2, \ldots, x_n) < 0 \quad (18)$$

and $v_{j_i^+}(t_{j \to j^+}) \leq \mu_{j \to j^+} v_{j_i}(t_{j \to j^+})$ (19)

where $\mu_{j \to j^+}$ are positive scalars.

Let us focus on the inequalities (19). Their aim is to ensure the global behavior of the like-Lyapunov function (17) at the switching time $t_{j \to j^+}$. These inequalities are verified if:

For $i = 1, \ldots, n$, $j_i = 1, \ldots, m_i$, $j_i^+ = 1, \ldots, m_i$, $s_{j_i} = 1, \ldots, r_{j_i}$:

$$E(X_{hj^+})^{-1} - \mu_{j \to j^+} E(X_{hj})^{-1} \leq 0 \quad (20)$$

That is to say:

$$(X^1_{hj^+})^{-1} - \mu_{j \to j^+} (X^1_{hj})^{-1} \leq 0 \quad (21)$$

Left and right multiplying by $X^1_{hj}$, and then using Schur complement, (21) is equivalent to

$$\begin{bmatrix} -\mu_{j \to j^+} X^1_{hj} & X^1_{hj} \\ X^1_{hj} & -X^1_{hj^+} \end{bmatrix} \leq 0 \quad (22)$$

Now, let us deal with (18), with the above defined notations, it can be rewritten as, $\forall t \neq t_{j \to j^+}$:

$$\sum_{i=1}^{n} \left[ \dot{\tilde{x}}_i^T E(X_{hj})^{-1} \tilde{x}_i + \tilde{x}_i^T (X_{hj})^{-1} E \dot{\tilde{x}}_i + \tilde{x}_i^T E(\dot{X}_{hj})^{-1} \tilde{x}_i \right] < 0 \quad (23)$$

Substituting (10) into (23), one can write, $\forall t \neq t_{j \to j^+}$:

$$\sum_{i=1}^{n} \begin{bmatrix} \tilde{x}_i^T \left[ \tilde{A}_{hj,hj}^T (X_{hj})^{-1} + (X_{hj})^{-1} \tilde{A}_{hj,hj} + E(\dot{X}_{hj})^{-1} \right] \tilde{x}_i \\ + \sum_{\alpha=1, \alpha \neq i}^{n} \tilde{x}_\alpha^T F_{i,\alpha,hj}^T (X_{hj})^{-1} \tilde{x}_i + \tilde{x}_i^T (X_{hj})^{-1} \tilde{F}_{i,\alpha,hj} \tilde{x}_\alpha \end{bmatrix} < 0 \quad (24)$$

From (12), the inequality (24) can be bounded by, $\forall t \neq t_{j \to j^+}$:

$$\sum_{i=1}^{n} \tilde{x}_i^T \begin{bmatrix} \tilde{A}_{hj,hj}^T (X_{hj})^{-1} + (X_{hj})^{-1} \tilde{A}_{hj,hj} + E(\dot{X}_{hj})^{-1} \\ + \sum_{\alpha=1, \alpha \neq i}^{n} \tau_{i,\alpha} (X_{hj})^{-1} \tilde{F}_{i,\alpha,hj} \tilde{F}_{i,\alpha,hj}^T (X_{hj})^{-1} \end{bmatrix} \tilde{x}_i + \sum_{\alpha=1, \alpha \neq i}^{n} \tau_{i,\alpha}^{-1} \tilde{x}_\alpha^T \tilde{x}_\alpha < 0 \quad (25)$$

Since $\sum_{i=1}^{n} \sum_{\alpha=1, \alpha \neq i}^{n} \tau_{i,\alpha}^{-1} x_\alpha^T x_\alpha = \sum_{i=1}^{n} \sum_{\alpha=1, \alpha \neq i}^{n} \tau_{\alpha,i}^{-1} x_i^T x_i$, $\forall x_i$, (25) is satisfied if, for $i = 1, \ldots, n$ and $\forall t \neq t_{j \to j^+}$:

$$\tilde{A}_{hj,hj}^T (X_{hj})^{-1} + (X_{hj})^{-1} \tilde{A}_{hj,hj} + E(\dot{X}_{hj})^{-1}$$
$$+ \sum_{\alpha=1, \alpha \neq i}^{n} \left[ \tau_{i,\alpha} (X_{hj})^{-1} \tilde{F}_{i,\alpha,hj} \tilde{F}_{i,\alpha,hj}^T (X_{hj})^{-1} + \tau_{\alpha,i}^{-1} I \right] < 0 \quad (26)$$

Such that $EX_{hj} = X_{hj} E > 0$; left and right multiplying the inequalities (26) respectively by $X_{hj}$, the equation (26) can be written as:

$$X_{hj} \tilde{A}_{hj,hj}^T + \tilde{A}_{hj,hj} X_{hj} + EX_{hj} (\dot{X}_{hj})^{-1} X_{hj} + \sum_{\alpha=1, \alpha \neq i}^{n} \left[ \tau_{i,\alpha} \tilde{F}_{i,\alpha,hj} \tilde{F}_{i,\alpha,hj}^T + \tau_{\alpha,i}^{-1} X_{hj} X_{hj} \right] < 0 \quad (27)$$

Now, the aim is to obtain the inequality (15) from (27). This can be achieved with the following usual mathematical developments. First, note that $-E(\dot{X}_{hj})^{-1} = E(X_{hj})^{-1} \dot{X}_{hj} (X_{hj})^{-1}$. This term is majored by $-\Phi_{s_{j_i} l_{j_i} k_{j_i} k'_{j_i}}$ [10]. Then to deals with the term $X_{hj} X_{hj}$, one apply the Schur complement.

***Condition 2** (Robustness condition)*: For all non-zero $\tilde{w}_{i,\alpha} \in L_2 [0 \; \infty)$, under zero initial condition $\tilde{x}_i(t_0) \equiv 0$, it holds that: for $i = 1, \ldots, N$,

$$\sum_{i=1}^{n} \left[ \dot{v}_i + \tilde{x}_i^T Q \tilde{x}_i^T - \varsigma_i^2 \sum_{\alpha=1, \alpha \neq i}^{N} \tilde{w}_{i,\alpha} \Xi \tilde{w}_{i,\alpha}^T \right] < 0 \quad (28)$$

$$\sum_{i=1}^{n} \begin{bmatrix} \dot{\tilde{x}}_i^T E(X_{hj})^{-1} \tilde{x}_i + \tilde{x}_i^T (X_{hj})^{-1} E \dot{\tilde{x}}_i + \tilde{x}_i^T E(\dot{X}_{hj})^{-1} \tilde{x}_i \\ + \tilde{x}_i^T Q \tilde{x}_i - \varsigma_i^2 \sum_{\alpha=1; \alpha \neq i}^{N} \tilde{w}_{i,\alpha}^T \Xi \tilde{w}_{i,\alpha} \end{bmatrix} < 0 \quad (29)$$

Substituting (10) into (29), one can write, $\forall t \neq t_{j \to j^+}$:

$$\sum_{i=1}^{n} \begin{bmatrix} \tilde{x}_i^T \left[ \tilde{A}_{hj,hj}^T (X_{hj})^{-1} + (X_{hj})^{-1} \tilde{A}_{hj,hj} + Q + E(\dot{X}_{hj})^{-1} \right] \tilde{x}_i \\ + \sum_{\alpha=1,\alpha \neq i}^{n} \begin{bmatrix} \tilde{x}_\alpha^T F_{i,\alpha,hj}^T (X_{hj})^{-1} \tilde{x}_i + \tilde{x}_i^T (X_{hj})^{-1} \tilde{F}_{i,\alpha,hj} \tilde{x}_\alpha - \varsigma_i^2 \tilde{w}_{i,\alpha}^T \Xi \tilde{w}_{i,\alpha} \\ + \tilde{w}_{i,\alpha}^T (\tilde{B}_{hj}^{w\alpha})^T (X_{hj})^{-1} \tilde{x}_i + \tilde{x}_i^T (X_{hj})^{-1} \tilde{B}_{hj}^{w\alpha} \tilde{w}_{i,\alpha} \end{bmatrix} \end{bmatrix} < 0 \quad (30)$$

From (12), the inequality (24) can be bounded by, $\forall t \neq t_{j \to j^+}$:

$$\sum_{i=1}^{n} \begin{bmatrix} \tilde{x}_i^T \begin{bmatrix} \tilde{A}_{hj,hj}^T (X_{hj})^{-1} + (X_{hj})^{-1} \tilde{A}_{hjhj} + Q + E(\dot{X}_{hj})^{-1} \\ + \sum_{\alpha=1,\alpha \neq i}^{n} \tau_{i,\alpha} (X_{hj})^{-1} \tilde{F}_{i,\alpha,hj} \tilde{F}_{i,\alpha,hj}^T (X_{hj})^{-1} \end{bmatrix} \tilde{x}_i \\ + \sum_{\alpha=1,\alpha \neq i}^{n} \left( \tilde{w}_{i,\alpha}^T (\tilde{B}_{hj}^{w\alpha})^T (X_{hj})^{-1} \tilde{x}_i + \tilde{x}_i^T (X_{hj})^{-1} \tilde{B}_{hj}^{w\alpha} \tilde{w}_{i,\alpha} \right) \\ + \sum_{\substack{\alpha=1 \\ \alpha \neq i}}^{n} \left( \tau_{i,\alpha}^{-1} \tilde{x}_\alpha^T \tilde{x}_\alpha - \varsigma_i^2 \tilde{w}_{i,\alpha}^T \Xi \tilde{w}_{i,\alpha} \right) \end{bmatrix} < 0 \quad (31):$$

Since $\sum_{i=1}^{n} \sum_{\alpha=1,\alpha \neq i}^{n} \tau_{i,\alpha}^{-1} x_\alpha^T x_\alpha = \sum_{i=1}^{n} \sum_{\alpha=1,\alpha \neq i}^{n} \tau_{\alpha,i}^{-1} x_i^T x_i$, $\forall x_i$ and $\forall t \neq t_{j \to j^+}$, (25) is satisfied if:

$$\sum_{i=1}^{n} \sum_{\alpha=1,\alpha \neq i}^{n} \begin{pmatrix} \tilde{x}_i^T \begin{pmatrix} N \left( \tilde{A}_{hjhj}^T (X_{hj})^{-1} + (X_{hj})^{-1} \tilde{A}_{hjhj} + Q + E(\dot{X}_{hj})^{-1} \right) \\ + \tau_{i,\alpha} (X_{h_{j_i}})^{-1} \tilde{F}_{i,\alpha,h_{j_i}} \tilde{F}_{i,\alpha,h_{j_i}}^T (X_{h_{j_i}})^{-1} + \tau_{\alpha,i}^{-1} I \end{pmatrix} \tilde{x}_i \\ + \tilde{w}_{i,\alpha}^T \tilde{B}_{w,\alpha,h_{j_i}}^T (X_{h_{j_i}})^{-1} \tilde{x}_i + \tilde{x}_i^T (X_{h_{j_i}})^{-1} \tilde{B}_{w,\alpha,h_{j_i}} \tilde{w}_{i,\alpha} \\ - \varsigma_i^2 \tilde{w}_{i,\alpha}^T \Xi \tilde{w}_{i,\alpha} \end{pmatrix} < 0 \quad (32)$$

The previous equation can be written as follow:

$$\begin{bmatrix} \tilde{x}_i \\ \tilde{w}_{i,\alpha} \end{bmatrix}^T \begin{bmatrix} \Upsilon_{hj,hj,hj} & * \\ (\tilde{B}_{hj}^{w\alpha})^T (X_{hj})^{-1} & -\varsigma_i^2 \Xi \end{bmatrix} \begin{bmatrix} \tilde{x}_i \\ \tilde{w}_{i,\alpha} \end{bmatrix} < 0 \quad (33)$$

with $\Upsilon_{hj,hj,hj} = \underline{N} \left( \tilde{A}_{hj,hj}^T (X_{hj})^{-1} + (X_{hj})^{-1} \tilde{A}_{hj,hj} + Q + E(\dot{X}_{hj})^{-1} \right)$
$+ \tau_{i,\alpha} (X_{hj})^{-1} \tilde{F}_{i,\alpha,hj} \tilde{F}_{i,\alpha,hj}^T (X_{hj})^{-1} + \tau_{\alpha,i}^{-1} I$

Left and right multiplying the inequalities (26) respectively by $\begin{bmatrix} X_{hj} & 0 \\ 0 & I \end{bmatrix}$, it yields for $i = 1,...,n$ and $\alpha = 1,...,n$ with $\alpha \neq i$:

$$\begin{bmatrix} \underline{N} \begin{pmatrix} X_{hj} \tilde{A}_{hjhj}^T + \tilde{A}_{hjhj} X_{hj} + X_{hj} Q X_{hj} \\ + E X_{hj} (\dot{X}_{hj})^{-1} X_{hj} \end{pmatrix} + \begin{pmatrix} \tau_{i,\alpha} \tilde{F}_{i,\alpha,hj} \tilde{F}_{i,\alpha,hj}^T \\ + \tau_{\alpha,i}^{-1} X_{hj} X_{hj} \end{pmatrix} & * \\ \hline X_{h_{j_i}} (\tilde{B}_{hj}^{w\alpha})^T & -\varsigma_i^2 \Xi \end{bmatrix} < 0 \quad (34)$$

To obtain the LMI condition (16), similarly to the previous part of this demonstration, from the property $-E(\dot{X}_{hj})^{-1} = E(X_{hj})^{-1} \dot{X}_{hj} (X_{hj})^{-1}$, we can major the derivative $-E\dot{X}_{hj}$ by $-\Phi_{s_{j_i} l_{j_i} k_{j_i} k'_{j_i}}$ [10] and then apply the Schur complement. ∎

## IV. NUMERICAL EXAMPLE

Let us consider the following system composed of two interconnected switched Takagi-sugeno subsystems given by:

**Subsystem 1**:

$$\dot{x}_1 = \sum_{j_1=1}^{2} \sum_{s_{j_1}=1}^{2} \xi_{j_1} h_{s_{j_1}} \left[ A_{s_{j_1}} x_1 + B_{s_{j_1}} u_1 + B_{s_{j_1}}^w w_1(t) + F_{1,2,s_{j_1}} x_2 + B_{s_{j_1}}^{w_2} w_2(t) \right] \quad (35)$$

$$y_1 = \sum_{j_1=1}^{2} \sum_{s_{j_1}=1}^{2} \xi_{j_1} h_{s_{j_1}} C_{s_{j_1}} x_1$$

with $x_1 = \begin{bmatrix} x_{11} \\ x_{12} \end{bmatrix}$ $A_{sj_1} = \begin{bmatrix} -2 & Ab_j \\ 0.1 & Aa_{sj} \end{bmatrix}$ $B_{sj_1} = \begin{bmatrix} Bb_j & Ba_{sj} \\ 0 & 1 \end{bmatrix}$

$C_{sj_1} = \begin{bmatrix} Ca_{sj} & 0.1 \\ -1 & 1 \end{bmatrix}$ $B_{sj_1}^w = \begin{bmatrix} wa_{sj} & wb_j \\ -.01 & .01 \end{bmatrix}$ $B_{sj_1}^{w_2} = \begin{bmatrix} .01 & \alpha b_j \\ \alpha a_{sj} & .01 \end{bmatrix}$

$F_{sj_1} = \begin{bmatrix} .01 & .01 & Fa_{sj} \\ Fb_j & .01 & .1 \end{bmatrix}$

In the **mode 1**, the variables values are given by $Ab_1 = 1$, $Aa_{11} = -2.1$, $Aa_{21} = -1.1$, $Bb_1 = -1.2$, $Ba_{11} = 0$, $Ba_{21} = 1.2$, $Ca_{11} = -.1$, $Ca_{12} = 1$, $Fb_1 = 0.01$, $Fa_{11} = .01$, $Fa_{21} = .1$, $wb_1 = 0.01$, $wa_{11} = -.01$, $wa_{21} = -.02$, $\alpha b_1 = 0.01$, $\alpha a_{11} = .02$, $\alpha a_{12} = .01$. In the **mode 2**, the variables values are given by: $Ab_2 = 0.2$, $Aa_{12} = -2$, $Aa_{22} = -3$, $Bb_2 = -1.5$, $Ba_{12} = 1$, $Ba_{22} = 3$, $Ca_{21} = 1$, $Ca_{22} = .1$, $Fb_2 = -0.01$, $Fa_{21} = .2$, $Fa_{22} = .02$, $wb_2 = -0.05$, $wa_{12} = -.05$, $wa_{22} = .01$, $\alpha b_2 = -0.05$, $\alpha a_{21} = .04$, $\alpha a_{22} = .03$. The membership functions are $h_{1_{1_1}}(x_1(t)) = \sin^2(x_{11}(t))$, $h_{2_{1_1}}(x_1(t)) = \sin^2(x_{12}(t))$, $h_{1_{2_1}}(x_1(t)) = 1 - h_{1_{1_1}}(x_1(t))$, and $h_{2_{2_1}}(x_1(t)) = 1 - h_{2_{1_1}}(x_1(t))$.

**Subsystem2**:

$$\dot{x}_2 = \sum_{j_2=1}^{2} \sum_{s_{j_2}=1}^{2} \xi_{j_2} h_{s_{j_2}} \left[ A_{s_{j_2}} x_2 + B_{s_{j_2}} u_2 + B_{s_{j_2}}^w w_2(t) + F_{2,1,s_{j_2}} x_1 + B_{s_{j_2}}^{w_1} w_1(t) \right] \quad (36)$$

$$y_2 = \sum_{j_2=1}^{2} \sum_{s_{j_2}=1}^{2} \xi_{j_2} h_{s_{j_2}} C_{s_{j_2}} x_2$$

with $x_2 = \begin{bmatrix} x_{21} \\ x_{22} \\ x_{23} \end{bmatrix}$ $A_{sj_2} = \begin{bmatrix} -2 & Ab_j & 0 \\ 0 & Aa_{sj} & 0 \\ 0 & .1 & -1.1 \end{bmatrix}$ $B_{sj_2} = \begin{bmatrix} -.1 & .5 & .1 \\ -.01 & .5 & .01 \\ -.01 & Ba_{sj} & .1 \end{bmatrix}$

$C_{sj_2} = \begin{bmatrix} .01 & Ca_{sj} & .1 \\ -1 & .1 & 1 \\ .1 & .1 & .1 \end{bmatrix}$ $F_{sj_2} = \begin{bmatrix} .01 & .001 & Fa_{sj} \\ .01 & .01 & Fb_j \end{bmatrix}$

$$B_{s_{j2}}^{w} = \begin{bmatrix} \alpha b_j & .05 & \alpha a_{sj} \\ .001 & .001 & .001 \\ .001 & .001 & .001 \end{bmatrix} \quad B_{s_{j2}}^{w_i} = \begin{bmatrix} wb_j & .05 & wa_{sj} \\ .001 & .001 & .001 \\ .001 & .001 & .001 \end{bmatrix}$$

In the **mode 1**, the variables values are given by $Ab_1 = 2$, $Aa_{11} = -1$, $Aa_{21} = -1.1$, $Ba_{11} = .01$, $Ba_{21} = .02$, $Ca_{11} = -.1$, $Ca_{12} = -.2$, $Fb_1 = 0.1$, $Fa_{11} = .2$, $Fa_{21} = .02$, $wb_1 = -0.01$, $wa_{11} = .01$, $wa_{21} = .001$, $\alpha b_1 = -.01$, $\alpha a_{11} = .01$, $\alpha a_{12} = .001$. In the **mode 2**, the variables values are given by: $Ab_2 = 1$, $Aa_{12} = -2$, $Aa_{22} = -3$, $Ba_{12} = 0.03$, $Ba_{22} = 0.04$, $Ca_{21} = -.4$, $Ca_{22} = -.3$, $Fb_2 = 0.2$, $Fa_{21} = Fa_{22} = .4$, $wb_2 = 0.01$, $wa_{12} = .02$, $wa_{22} = .05$, $\alpha b_2 = .01$, $\alpha a_{21} = .02$, $\alpha a_{22} = .05$ and the membership functions are $h_{1_{l_2}}(x_2) = \sin^2(x_{21})$, $h_{1_{2_2}}(x_2) = 1 - h_{1_{l_2}}(x_2)$ $h_{2_{l_2}}(x_2) = \sin^2(x_{22})$, and $h_{2_{2_2}}(x_2) = 1 - h_{2_{l_2}}(x_2)$.

Let us assume that each subsystem switches within the frontier defined by: $H_{11} = 0.9x_{11} + x_{12}$, $H_{12} = -0.2x_{11} + 9x_{12}$, $H_{21} = -x_{21} + x_{22}$ and $H_{22} = x_{21} - 2x_{22}$. The external disturbances $w_1$ and $w_2$ are considered as white noise sequences.

A set of decentralized switched controllers (5) is synthesized based on theorem 1 via the Matlab LMI toolbox. To do so, the lower bounds of membership functions are chosen as $\lambda_{1_{l_1}} = \lambda_{1_{2_1}} = \lambda_{1_{l_2}} = \lambda_{1_{2_2}} = -6$ and the disturbance attenuation levels are prefixed by $\varsigma_1^2 = 1.7$ and $\varsigma_2^2 = 1.5$.

The close-loop subsystem dynamics are shown in Figure 1 for the initial states $x_1(0) = [2 \ 2]^T$ and $x_2(0) = [-1 \ 1.5 \ -1]^T$.

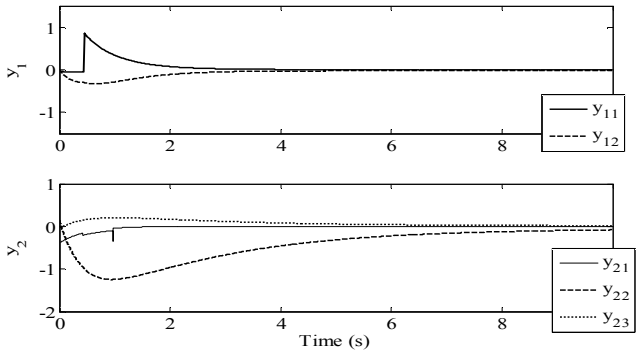

Figure 1. States dynamics of the overall closed-loop interconnected switched Takagi-Sugeno system.

Figure 2 shows the control signals as well as the switching modes' evolution. As expected, the synthesized decentralized switched controller stabilizes the overall large scale switched system composed of (1) and (2).

## V. CONCLUSION

This study has focused on large scale switched nonlinear systems where each nonlinear mode is represented by Takagi-Sugeno systems with external disturbances. To ensure the stability of the whole closed-loop system, a set of decentralized static output feedback switched non-PDC controllers has been considered. Therefore, LMI-based conditions for the design of the decentralized controller have been proposed through the consideration of a multiple switched non quadratic Lyapunov function and by using a descriptor redundancy formulation. Finally, a numerical example has been proposed to show the effectiveness of the proposed approach. An extension of the proposed approach to general switched systems under asynchronous switching will be the focus of our future work.

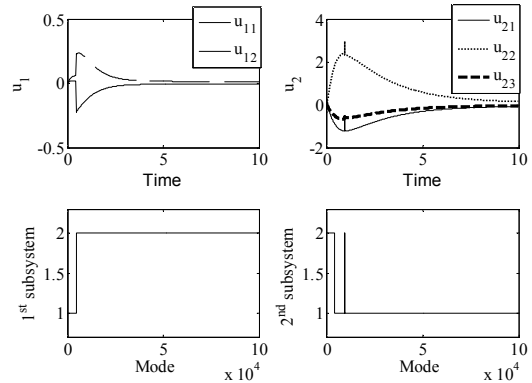

Figure 2. Control signals and switched laws' evolutions.